\documentstyle[amssymb,12pt]{article}
\newtheorem{th}{Theorem}[section]

\newtheorem{cor}[th]{Corollary}
\newtheorem{defn}[th]{Definition}
\newenvironment{defn-new}{\begin{defn} \em}{\end{defn}}
\newtheorem{rem}[th]{Remark}
\newenvironment{rem-new}{\begin{rem} \em}{\end{rem}}
\newtheorem{ex}[th]{Example}
\newenvironment{ex-new}{\begin{ex} \em}{\end{ex}}

\makeatletter \@addtoreset{equation}{section} \makeatother
\begin{document}

\setcounter{page}{1}
\thispagestyle{empty}

\bigskip

\begin{center}
{\Large {\bf Ricci solitons in contact metric manifolds}}\bigskip \bigskip

{Mukut Mani Tripathi}\bigskip

Department of Mathematics

Banaras Hindu University

Varanasi 221 005, India.

Email: {\tt mmtripathi66@yahoo.com}
\end{center}

\bigskip

\begin{quote}
\noindent {\bf Abstract.} In $N(k)$-contact metric manifolds and/or $\left(
k,\mu \right) $-manifolds, gradient Ricci solitons, compact Ricci solitons
and Ricci solitons with $V$ pointwise collinear with the structure vector
field $\xi $ are studied.


\end{quote}

\section{Introduction\label{sect-intro}}

A Ricci soliton is a generalization of an Einstein metric. In a Riemannian
manifold $\left( M,g\right) $, $g$ is called a Ricci soliton \cite%
{Hamilton-88} if
\begin{equation}
{\pounds }_{V}g+2\,Ric+2\lambda g=0,  \label{eq-Ricci-soliton}
\end{equation}%
where $\pounds $ is the Lie derivative, $V$ is a complete vector field on $M$
and $\lambda $ is a constant. Metrics satisfying (\ref{eq-Ricci-soliton})
are interesting and useful in physics and are often referred as
quasi-Einstein (e.g. \cite{Chave-Valent-1}, \cite{Chave-Valent-2}, \cite%
{Friedan-85}). Compact Ricci solitons are the fixed point of the Ricci flow
\[
\frac{\partial }{\partial t}\,g=-\,2\,Ric
\]%
projected from the space of metrics onto its quotient modulo diffeomorphisms
and scalings, and often arise as blow-up limits for the Ricci flow on
compact manifolds. The Ricci soliton is said to be shrinking, steady, and
expanding according as $\lambda $ is negative, zero, and positive
respectively. If the vector field $V$ is the gradient of a potential
function $-f$, then $g$ is called a gradient Ricci soliton and equation (\ref%
{eq-Ricci-soliton}) assumes the form
\begin{equation}
{\frak \nabla \nabla }f=Ric+\lambda g.  \label{eq-Ricci-soliton-grad}
\end{equation}%
A Ricci soliton on a compact manifold has constant curvature in dimension $2$
(Hamilton \cite{Hamilton-88}), and also in dimension $3$ (Ivey \cite{Ivey-93}%
). For details we refer to Chow and Knoff \cite{Chow-Knoff} and Derdzinski
\cite{Derdzinski}. We also recall the following significant result of
Perelman \cite{Perelman}: {\em A Ricci soliton on a compact manifold is a
gradient Ricci soliton}.

\medskip On the other hand, the roots of contact geometry lie in
differential equations as in 1872 Sophus Lie introduced the notion of
contact transformation (Ber\"{u}hrungstransformation) as a geometric tool to
study systems of differential equations. This subject has manifold
connections with the other fields of pure mathematics, and substantial
applications in applied areas such as mechanics, optics, phase space of a
dynamical system, thermodynamics and control theory (for more details see
\cite{Arnold-89}, \cite{Blair-02}, \cite{Geiges-01}, \cite{MacLane} and \cite%
{Naz-Shat-Ster}).

\medskip It is well known \cite{Sasaki-62} that the tangent sphere bundle $%
T_{1}M$ of a Riemannian manifold $M$ admits a contact metric structure. If $%
M $ is of constant curvature $c=1$ then $T_{1}M$ is Sasakian \cite%
{Tashiro-69}, and if $c=0$\ then the curvature tensor $R$ satisfies $%
R(X,Y)\xi =0\,$\ \cite{Blair-77}. As a generalization of these two cases, in
\cite{BKP-95}, Blair, Koufogiorgos and Papantoniou started the study of the
class of contact metric manifolds, in which the structure vector field $\xi$
satisfies the $\left( k,\mu \right) $-nullity condition. A contact metric
manifold belonging to this class is called a $(k,\mu)$-manifold. Such a
structure was first obtained by Koufogiorgos \cite{Koufog-93} by applying a $%
D_{a}$-homothetic deformation \cite{Tanno-68} on a contact metric manifold
satisfying $R(X,Y)\xi =0$. In particular, a $(k,0)$-manifold is called an $%
N(k)$-contact metric manifold (\cite{BBK-92}, \cite{BKT-05}, \cite{Tanno-88}%
) and generalizes the cases $R(X,Y)\xi =0$, $K$-contact and Sasakian.

\medskip In \cite{Sharma-07}, Sharma has started the study of Ricci solitons
in $K$-contact manifolds. In a $K$-contact manifold the structure vector
field $\xi $ is Killing, that is, ${\pounds }_{\xi }g=0$; which is not in
general true in contact metric manifolds. Motivated by these circumstances,
in this paper we study Ricci solitons in $N(k)$-contact metric manifolds and
$\left( k,\mu \right) $-manifolds. In section~\ref{sect-contact-mfd}, we
give a brief description of $N(k)$-contact metric manifolds and $\left(
k,\mu \right) $-manifolds. In section~\ref{sect-main-results}, we prove main
results.

\section{Contact metric manifolds \label{sect-contact-mfd}}

A $1$-form $\eta $ on a $\left( 2n+1\right) $-dimensional smooth manifold $M$
is called a {\em contact form} if $\eta \wedge (d\eta )^{n}\neq 0$
everywhere on $M$, and $M$ equipped with a contact form is a {\em contact
manifold}. For a given contact $1$-form $\eta $, there exists a unique
vector field $\xi $, called the {\em characteristic vector field}, such that
$\eta (\xi )=1$, $d\eta (\xi ,\cdot )=0$, and consequently ${\pounds }_{\xi
}\eta =0$, ${\pounds }_{\xi }d\eta =0$. In 1953, Chern \cite{Chern-53}
proved that the structural group of a $\left( 2n+1\right) $-dimensional
contact manifold can be reduced to ${\cal U}\left( n\right) \times 1$. A $%
(2n+1)$-dimensional differentiable manifold $M$ is called an {\em almost
contact manifold} \cite{Gray-59} if its structural group can be reduced to $%
{\cal U}\left( n\right) \times 1$. Equivalently, there is an {\em almost
contact structure} $\left( \varphi ,\xi ,\eta \right) $ \cite{Sasaki-60}
consisting of a tensor field $\varphi $ of type $\left( 1,1\right) $, a
vector field $\xi $, and a $1$-form $\eta $ satisfying
\begin{equation}
\varphi ^{2}=-I+\eta \otimes \xi ,\quad \eta (\xi )=1,\quad \varphi \xi
=0,\quad \eta \circ \varphi =0.  \label{eq-phi-eta-xi}
\end{equation}%
First and one of the remaining three relations of (\ref{eq-phi-eta-xi})
imply the other two relations. An almost contact structure is {\em normal}
\cite{SH-61} if the torsion tensor $\left[ \varphi ,\varphi \right] +2d\eta
\otimes \xi $, where $\left[ \varphi ,\varphi \right] $ is the Nijenhuis
tensor of $\varphi $, vanishes identically. Let $g$ be a compatible
Riemannian metric with $\left( \varphi ,\xi ,\eta \right) $, that is,
\begin{equation}
g\left( X,Y\right) =g\left( \varphi X,\varphi Y\right) +\eta \left( X\right)
\eta \left( Y\right) ,\qquad X,Y\in TM.  \label{eq-metric-1}
\end{equation}%
Then, $M$ becomes an {\em almost contact metric manifold} equipped with an
{\em almost contact metric structure} $\left( \varphi ,\xi ,\eta ,g\right) $%
. The equation (\ref{eq-metric-1}) is equivalent to
\begin{equation}
g\left( X,\varphi Y\right) =-g\left( \varphi X,Y\right) \quad {\rm alongwith}%
\quad g\left( X,\xi \right) =\eta \left( X\right) .  \label{eq-metric-2}
\end{equation}%
An almost contact metric structure becomes a contact metric structure if $%
g\left( X,\varphi Y\right) =d\eta (X,Y)$ for all $X,Y\in TM$. In a contact
metric manifold $M$, the $\left( 1,1\right) $-tensor field $h$ defined by $%
2h={\pounds }_{\xi }\varphi $, is symmetric and satisfies
\begin{equation}
h\xi =0,\qquad h\varphi +\varphi h=0,  \label{eq-cont-h}
\end{equation}%
\begin{equation}
\nabla \xi =-\ \varphi -\varphi h,  \label{eq-cont-del-xi}
\end{equation}%
where $\nabla $ is the Levi-Civita connection. A contact metric manifold is
called a $K$-{\em contact manifold\/} if the characteristic vector field $%
\xi $ is a Killing vector field. An almost contact metric manifold is a $K$%
{\em -contact manifold} if and only if $\nabla \xi =-\varphi $. A $K$%
-contact manifold is a contact metric manifold, while the converse is true
if $h=0$. A normal contact metric manifold is a {\em Sasakian manifold}. A
contact metric manifold $M$ is Sasakian if and only if the curvature tensor $%
R$ satisfies
\begin{equation}
R(X,Y)\xi =\eta (Y)X-\eta (X)Y,\qquad X,Y\in TM.  \label{eq-Sas-2}
\end{equation}%
A contact metric manifold $M$ is said to be $\eta ${\em -Einstein} (\cite%
{Okumura-62} or see \cite{Blair-02} p. 105) if the Ricci tensor $Ric$
satisfies $Ric=ag+b\eta \otimes \eta $, where $a$ and $b$ are some smooth
functions on the manifold. In particular if $b=0$, then $\,M\,$ becomes an
{\em Einstein manifold}.

\medskip A Sasakian manifold is always a $K$-contact manifold. The converse
is true if either the dimension is three (\cite{Blair-02}, p.~76), or it is
compact Einstein (Theorem A, \cite{Boyer-Galicki}) or compact $\eta $%
-Einstein with $a>-2$ (Theorem 7.2, \cite{Boyer-Galicki}). The conclusions
of Theorems A and 7.2 of \cite{Boyer-Galicki} are still true if the
condition of compactness is weakened to completeness (Proposition 1, \cite%
{Sharma-07}).

\medskip In \cite{BKP-95}, Blair, Koufogiorgos and Papantoniou introduced a
class of contact metric manifolds $M$, which satisfy
\begin{equation}
R(X,Y)\xi =\left( kI+\mu h\right) \left( \eta \left( Y\right) X-\eta \left(
X\right) Y\right) ,\qquad X,Y\in TM,  \label{eq-km}
\end{equation}%
where $k$, $\mu $ are real constants. A contact metric manifold belonging to
this class is called a $\left( k,\mu \right) ${\em -manifold}. If $\mu =0$,
then a $(k,\mu )$-manifold is called an $N(k)$-{\em contact metric manifold}
(\cite{BBK-92}, \cite{BKT-05}, \cite{Tanno-88}). In a $(k,\mu )$-manifold $M$%
, one has \cite{BKP-95} 
\begin{eqnarray}
(\nabla _{X}h)Y &=&\left( \left( 1-k\right) g\left( X,\varphi Y\right)
+g\left( X,\varphi hY\right) \right) \xi  \nonumber \\
&&+\,\eta \left( Y\right) \left( h\left( \varphi X+\varphi hX\right) \right)
-\mu \eta \left( X\right) \varphi hY  \label{eq-derivative-h}
\end{eqnarray}%
for all $X,Y\in TM$. The Ricci operator $Q$ satisfies $Q\xi =2nk\xi $, where
$\dim (M)=2n+1$. Moreover, $h^{2}=\left( k-1\right) \varphi ^{2}$ and $k\leq
1$. In fact, for a $(k,\mu )$-manifold, the conditions of being a Sasakian
manifold, a $K$-contact manifold, $k=1$ and $h=0$ are all equivalent. The
tangent sphere bundle $T_{1}M$ of a Riemannian manifold $M$ of constant
curvature $c$ is a $(k,\mu )$-manifold with $k=c(2-c)$ and $\mu =-2c$.
Characteristic examples of non-Sasakian $(k,\mu )$-manifolds are the tangent
sphere bundles of Riemannian manifolds of constant curvature not equal to
one and certain Lie groups \cite{Boeckx-00}. For more details we refer to
\cite{Blair-02} and \cite{BKP-95}.

\section{Main results\label{sect-main-results}}

Let $(M,\varphi ,\xi ,\eta ,g)$ be a $\left( 2n+1\right) $-dimensional
non-Sasakian $\left( k,\mu \right) $-manifold. Then the Ricci operator $Q$
is given by \cite{BKP-95}
\begin{equation}
Q=2nkI+\left( 2\left( n-1\right) +\mu \right) h-\left( 2\left( n-1\right)
-n\mu +2nk\right) \varphi ^{2}.  \label{eq-Q-1}
\end{equation}%
We also have
\begin{equation}
\left( \nabla _{X}\varphi ^{2}\right) Y=\left( X\eta \left( Y\right) \right)
\xi -\eta \left( \nabla _{X}Y\right) \xi -\eta \left( Y\right) \varphi
X-\eta \left( Y\right) \varphi hX,  \label{eq-der-phi-sq}
\end{equation}%
where first equation of (\ref{eq-phi-eta-xi}) and equation (\ref%
{eq-cont-del-xi}) are used. Using (\ref{eq-derivative-h}) and (\ref%
{eq-der-phi-sq}) from (\ref{eq-Q-1}) we obtain
\begin{eqnarray*}
\left( \nabla _{X}Q\right) Y &=&\left( 2\left( n-1\right) +\mu \right)
\left\{ \left( 1-k\right) g\left( X,\varphi Y\right) \xi +g\left( X,\varphi
hY\right) \xi -\mu \eta \left( X\right) \varphi hY\right\} \\
&&-\left( 2(n-1)-n\mu +2nk\right) \left\{ \left( X\eta \left( Y\right)
\right) \xi -\eta \left( \nabla _{X}Y\right) \xi \right\} \\
&&+\left( 2\left( 2n-1\right) k-\left( n+1\right) \mu +k\mu \right) \eta
\left( Y\right) \varphi X+\left( \left( n+1\right) \mu -2nk\right) \eta
\left( Y\right) h\varphi X.
\end{eqnarray*}%
Consequently, we have
\begin{eqnarray}
\left( \nabla _{X}Q\right) Y-\left( \nabla _{Y}Q\right) X &=&\left( 2\left(
n+1\right) \mu -4\left( 2n-1\right) k-2k\mu \right) d\eta \left( X,Y\right)
\xi  \nonumber \\
&&+\left( 2\left( 2n-1\right) k-\left( n+1\right) \mu +k\mu \right) \left(
\eta \left( Y\right) \varphi X-\eta \left( X\right) \varphi Y\right)
\nonumber \\
&&+\left( \left( \mu +3n-1\right) \mu -2nk\right) \left( \eta \left(
Y\right) \varphi hX-\eta \left( X\right) \varphi hY\right) ,
\label{eq-Q-der}
\end{eqnarray}%
where (\ref{eq-metric-2}) has been used.

\medskip We also recall the following results for later use.

\begin{th}
\label{th-Tanno-88} {\rm (Theorem~5.2, Tanno \cite{Tanno-88})} An Einstein $%
N(k )$-contact metric manifold of dimension $\ge 5$ is necessarily Sasakian.
\end{th}

\begin{th}
\label{th-non-Sas-Einstein} {\rm (Theorem 1.2, Tripathi and Kim \cite%
{Tri-Kim-04})} A non-Sasakian Einstein $(k,\mu )$-manifold is flat and $3$%
-dimensional.
\end{th}

\medskip Now we prove the following

\begin{th}
\label{th-grad-soliton} If the metric $g$ of an $N(k)$-contact metric
manifold $(M,g)$ is a gradient Ricci soliton, then

\begin{enumerate}
\item[{\bf (a)}] either the potential vector field
is a nullity vector field,

\item[{\bf (b)}] or $g$ is a shrinking soliton and $(M,g)$ is Einstein
Sasakian,

\item[{\bf (c)}] or $g$ is a steady soliton and $(M,g)$ is $\,3$-dimensional
and flat.
\end{enumerate}
\end{th}

\noindent {\bf Proof.} Let $(M,g)$ be a $\left( 2n+1\right) $-dimensional $%
N(k)$-contact metric manifold and $g$ a gradient Ricci soliton. Then the
equation (\ref{eq-Ricci-soliton-grad}) can be written as
\begin{equation}
\nabla _{Y}Df=QY+\lambda Y  \label{eq-Ricci-soliton-grad-2}
\end{equation}%
for all vector fields $Y$ in $M$, where $D$ denotes the gradient operator of
$g$. From (\ref{eq-Ricci-soliton-grad-2}) it follows that
\begin{equation}
R\left( X,Y\right) Df=\left( \nabla _{X}Q\right) Y-\left( \nabla
_{Y}Q\right) X,\qquad X,Y\in TM.  \label{eq-R(X,Y)Df}
\end{equation}%
We have
\begin{equation}
g\left( R\left( \xi ,Y\right) Df,\xi \right) =g\left( k\left( Df-\left( \xi
f\right) \xi \right) ,Y\right) ,\qquad Y\in TM,  \label{eq-R(xi,Y)Df}
\end{equation}%
where (\ref{eq-km}) with $\mu =0$ is used. Also in an $N\left( k\right) $%
-contact metric manifold, it follows that
\begin{equation}
g\left( \left( \nabla _{\xi }Q\right) Y-\left( \nabla _{Y}Q\right) \xi ,\xi
\right) =0,\qquad Y\in TM.  \label{eq-Q-der-xi}
\end{equation}%
From (\ref{eq-R(X,Y)Df}), (\ref{eq-R(xi,Y)Df}) and (\ref{eq-Q-der-xi}) we
get
\[
k\left( Df-\left( \xi f\right) \xi \right) =0,
\]%
that is, either $k=0$ or
\begin{equation}
Df=\left( \xi f\right) \xi .  \label{eq-Df}
\end{equation}%

If $k=0$, then putting $k=0=\mu $ in (\ref{eq-Q-der}), it follows that $Q$
is a Codazzi tensor, that is,
\[
\left( \nabla _{X}Q\right) Y-\left( \nabla _{Y}Q\right) X=0,\qquad X,Y\in
TM,
\]%
which in view of (\ref{eq-R(X,Y)Df}) gives
\[
R\left( X,Y\right) Df=0,\qquad X,Y\in TM,
\]%
that is, the potential vector field $\,{Df}\,$ is a nullity vector field
(see \cite{Clift-Malt} and \cite{Tanno-78} for details).

\medskip Now, we assume that (\ref{eq-Df}) is true. Using (\ref{eq-Df}) in (%
\ref{eq-Ricci-soliton-grad-2}) we get
\[
Ric\left( X,Y\right) +\lambda g\left( X,Y\right) =Y\left( \xi f\right) \eta
\left( X\right) -\left( \xi f\right) g\left( X,\varphi Y\right) -\left( \xi
f\right) g\left( X,\varphi hY\right) ,
\]%
where (\ref{eq-cont-del-xi}) is used. Symmetrizing this with respect to $X$
and $Y$ we obtain
\begin{equation}
2\,Ric\left( X,Y\right) +2\lambda g\left( X,Y\right) =X\left( \xi f\right)
\eta \left( Y\right) +Y\left( \xi f\right) \eta \left( X\right) -2\left( \xi
f\right) g\left( \varphi hX,Y\right) .  \label{eq-2Ric}
\end{equation}%
Putting $Y=\xi $, we get
\begin{equation}
X\left( \xi f\right) =\left( 2nk+\lambda \right) \eta \left( X\right) .
\label{eq-X(xi-f)}
\end{equation}%
From (\ref{eq-2Ric}) and (\ref{eq-X(xi-f)}) we get
\begin{equation}
Ric\left( X,Y\right) +\lambda g\left( X,Y\right) =\left( 2nk+\lambda \right)
\eta \left( X\right) \eta \left( Y\right) -\left( \xi f\right) g\left(
\varphi hX,Y\right) .  \label{eq-Ric}
\end{equation}%
Using (\ref{eq-Ric}) in (\ref{eq-Ricci-soliton-grad-2}), we get
\begin{equation}
\nabla _{Y}Df=\left( 2nk+\lambda \right) \eta \left( Y\right) \xi -\left(
\xi f\right) \varphi hY.  \label{eq-Ricci-soliton-grad-3}
\end{equation}%
Using (\ref{eq-Ricci-soliton-grad-3}) we compute $R\left( X,Y\right) Df$ and
obtain
\begin{equation}
g\left( R\left( X,Y\right) \left( \xi f\right) \xi ,\xi \right) =4\left(
2nk+\lambda \right) d\eta \left( X,Y\right) ,  \label{eq-R(X,Y)(xi-f)xi}
\end{equation}%
where equations (\ref{eq-Df}) and (\ref{eq-cont-del-xi}) are used. Thus we
get
\begin{equation}
2nk+\lambda =0  \label{eq-2nk-plus-lambda}
\end{equation}%
Therefore from equation (\ref{eq-X(xi-f)}) we have
\[
X\left( \xi f\right) =0,\qquad X\in TM,
\]%
that is,
\[
\xi f=c,
\]%
where $c$ is a constant. Thus the equation (\ref{eq-Df}) gives
\[
df=c\,\eta \,.
\]%
Its exterior derivative implies that
\[
c\,d\eta =0,
\]%
that is, $\,c=0$. Hence $f$ is constant. Consequently, the equation (\ref%
{eq-Ricci-soliton-grad-2}) reduces to
\[
Ric=-\lambda g=2nkg,
\]%
that is, $M$ is Einstein. Then in view of Theorem~\ref{th-non-Sas-Einstein}
and Theorem~\ref{th-Tanno-88}, it follows that either $M$ is Sasakian or $M$
is $3$-dimensional and flat. In case of Sasakian, $\lambda =-2n$ is
negative, and therefore the soliton $g$ is shrinking. In case of $3$%
-dimensional and flat, $\lambda =0$, and therefore the soliton $g$ is
steady. $\blacksquare $

\begin{cor}
\label{cor-grad-soliton} Let $(M,g)$ be a compact $N(k)$-contact metric
manifold with $k\neq 0$. If $g$ is a Ricci soliton, then $g$ is a shrinking
soliton and $(M,g)$ is Einstein Sasakian.
\end{cor}

\noindent {\bf Proof.} The proof follows from Theorem~\ref{th-grad-soliton}
and the following significant result of Perelman \cite{Perelman}: A Ricci
soliton on a compact manifold is a gradient Ricci soliton. $\blacksquare $

\medskip In \cite{Sharma-07}, a corollary of Theorem 1 is stated as follows:
If the metric $g$ of a compact $K$-contact manifold is a Ricci soliton, then
$g$ is a shrinking soliton which is Einstein Sasakian. In Corollary~\ref%
{cor-grad-soliton}, the assumptions are weakened.

\medskip Next, we have the following

\begin{th}
In a non-Sasakian 
$\,\left( k,\mu \right) \,$-manifold $(M,g)$ if $g$ is a compact Ricci
soliton, then $(M,g)$ is $3$-dimensional and flat.
\end{th}

\noindent {\bf Proof.} In a non-Sasakian $\left( k,\mu \right) $-manifold,
the scalar curvature $r$ is given by \cite{BKP-95}
\begin{equation}
r=2n\left( 2n-2+k-n\mu \right) .  \label{eq-scalar}
\end{equation}%
Consequently, the scalar curvature is a constant. If $g$ is a compact Ricci
soliton, then by Proposition~2 of \cite{Sharma-07}, which states that a
compact Ricci soliton of constant scalar curvature is Einstein, it follows
that the non-Sasakian $\left( k,\mu \right) $-manifold is Einstein. Then by
Theorem~\ref{th-non-Sas-Einstein}, it becomes $3$-dimensional and flat. $%
\blacksquare $

\medskip Given a non-Sasakian $(\kappa ,\mu )$-manifold $M$, Boeckx \cite%
{Boeckx-00} introduced an invariant
\[
I_{M}=\frac{1-\mu /2}{\sqrt{1-\kappa }}
\]%
and showed that for two non-Sasakian $(\kappa ,\mu )$-manifolds $%
(M_{i},\varphi _{i},\xi _{i},\eta _{i},g_{i})$, $i=1,2$, we have $%
I_{M_{1}}=I_{M_{2}}$ if and only if up to a $D$-homothetic deformation, the
two manifolds are locally isometric as contact metric manifolds. Thus we
know all non-Sasakian $(\kappa ,\mu )$-manifolds locally as soon as we have
for every odd dimension $2n+1$ and for every possible value of the invariant
$I$, one $(\kappa ,\mu )$-manifold $(M,\varphi ,\xi ,\eta ,g)$ with $I_{M}=I$%
. For $I>-1$ such examples may be found from the standard contact metric
structure on the tangent sphere bundle of a manifold of constant curvature $%
c $ where we have $I=\frac{1+c}{|1-c|}$. Boeckx also gives a Lie algebra
construction for any odd dimension and value of $I\leq -1$.

\medskip In the following, we recall Example~3.1 of \cite{BKT-05}.

\begin{ex-new}
\label{example-BKT} For $n>1$, the Boeckx invariant for a $(2n+1)$%
-dimensional $\left( 1-\frac{1}{n},0\right) $-manifold is $\sqrt{n}>-1$.
Therefore, we consider the tangent sphere bundle of an $(n+1)$-dimensional
manifold of constant curvature $c$ so chosen that the resulting $D_{a}$%
-homothetic deformation will be a $\left( 1-\frac{1}{n},0\right) $-manifold.
That is for $k =c(2-c)$ and $\mu =-2c$ we solve
\[
1-\frac{1}{n}=\frac{k +a^{2}-1}{a^{2}},\quad 0=\frac{\mu +2a-2}{a}
\]%
for $a$ and $c$. The result is
\[
c=\frac{\left( \sqrt{n}\pm 1\right) ^{2}}{n-1},\quad a=1+c
\]%
and taking $c$ and $a$ to be these values we obtain a $N\left( 1-\frac{1}{n}%
\right) $-contact metric manifold.
\end{ex-new}

In \cite{Sharma-07}, Sharma noted that if a $K$-contact metric is a Ricci
soliton with $V=\xi $ then it is Einstein. Even in more general case, he
showed that if a $K$-contact metric is a Ricci soliton with $V$ pointwise
collinear with $\xi $ then $V$ is a constant multiple of $\xi $ (hence
Killing) and $g$ is Einstein. Here we prove the following

\begin{th}
Let $\,(M,g)\,$ be a non-Sasakian $($or non-$K$-contact\/$)$ $\,N\left(
k\right) $-contact metric manifold. If the metric $\,g\,$ is a Ricci soliton
with $V$ pointwise collinear with $\xi $, then $\dim (M)>3$, the metric $g$
is a shrinking Ricci soliton and $M$ is locally isometric to a contact
metric manifold obtained by a $D_{\left( 1+\frac{\left( \sqrt{n}\pm 1\right)
^{2}}{n-1}\right) }$-homothetic deformation of the contact metric structure
on the tangent sphere bundle of an $\left( n+1\right) $-dimensional
Riemannian manifold of constant curvature $\frac{\left( \sqrt{n}\pm 1\right)
^{2}}{n-1}$.
\end{th}

\noindent {\bf Proof.} Let $(M,g)$ be a $\left( 2n+1\right) $-dimensional
contact metric manifold and the metric $g$ a Ricci soliton with $V=\alpha
\xi $ ($\alpha $ being a function on $M$). Then from (\ref{eq-Ricci-soliton}%
) we obtain
\begin{equation}
2Ric\left( X,Y\right) =-2\lambda g\left( X,Y\right) +2\alpha g\left( \varphi
hX,Y\right) -g\left( \left( X\alpha \right) \xi ,Y\right) -g\left( X,\left(
Y\alpha \right) \xi \right) ,  \label{eq-Ric-01}
\end{equation}%
where (\ref{eq-cont-del-xi}) and (\ref{eq-metric-2}) are used. Now let $(M,g)
$ be an $N\left( k\right) $-contact metric manifold. Putting $X=\xi =Y$ in (%
\ref{eq-Ric-01}) and using $h\xi =0$ and $Q\xi =2nk$ we get
\begin{equation}
\xi \alpha +2nk+\lambda =0.  \label{eq-xi-alpha}
\end{equation}%
Again putting $X=\xi $ in (\ref{eq-Ric-01}) and using $h\xi =0$, $Q\xi =2nk$
and (\ref{eq-xi-alpha}) we get
\begin{equation}
d\alpha =\left( 2nk+\lambda \right) \eta ,  \label{eq-d-alpha}
\end{equation}%
which shows that $\alpha $ is a constant and $\lambda =-2nk$; and
consequently (\ref{eq-Ric-01}) becomes
\begin{equation}
Ric\left( X,Y\right) =2nkg\left( X,Y\right) +\alpha g\left( \varphi
hX,Y\right) .  \label{eq-Ric-02}
\end{equation}%
At this point, we assume that $(M,g)$ is also non-Sasakian. It is known that
in a $(2n+1)$-dimensional non-Sasakian $(k,\mu )$-manifold $M$ the Ricci
tensor is given by \cite{BKP-95}
\begin{eqnarray}
Ric\left( X,Y\right)  &=&\left( 2\left( n-1\right) -n\mu \right) g\left(
X,Y\right) +\left( 2\left( n-1\right) +\mu \right) g\left( hX,Y\right)
\nonumber \\
&&+\ \left( 2\left( 1-n\right) +n\left( 2k+\mu \right) \right) \eta \left(
X\right) \eta \left( Y\right) .  \label{eq-Ric-curvature}
\end{eqnarray}%
Consequently, putting $\mu =0$ in (\ref{eq-Ric-curvature}) we get
\begin{eqnarray}
Ric\left( X,Y\right)  &=&2\left( n-1\right) g\left( X,Y\right) +2\left(
n-1\right) g\left( hX,Y\right)   \nonumber \\
&&+\left( 2\left( 1-n\right) +2nk\right) \eta \left( X\right) \eta \left(
Y\right) .  \label{eq-Ric-03}
\end{eqnarray}%
Replacing $X$ by $\varphi X$ in equations (\ref{eq-Ric-02}) and (\ref%
{eq-Ric-03}) and equating the right hand sides of the resulting equations we
get
\begin{equation}
\left( 2nk-2\left( n-1\right) \right) g\left( \varphi X,Y\right) =\alpha
g\left( hX,Y\right) +2\left( n-1\right) g\left( \varphi hX,Y\right) .
\label{eq-Ric-04}
\end{equation}%
If $n=1$, from (\ref{eq-Ric-04}) we get
\[
2kg\left( \varphi X,Y\right) =\alpha g\left( hX,Y\right) ,
\]%
which gives $h=0$, a contradiction. If $n>1$, anti-symmetrizing the equation
(\ref{eq-Ric-04}) we get
\[
nk-n+1=0,
\]%
which gives $k=1-1/n$. Using $n>1$ and $k=1-1/n$ in $\lambda =-2nk$, we get $%
\lambda =2\left( 1-n\right) <0$, which shows that $g$ is a shrinking Ricci
soliton. Finally, in view of $n>1$, $k=1-1/n$ and the Example~\ref%
{example-BKT}, the proof is complete. $\blacksquare $

\bigskip \noindent {\bf Acknowledgement:} The author is thankful to
Professor Ramesh Sharma, University of New Haven, USA for some
useful discussion during the preparation of this paper.

\end{document}